\documentstyle{article}[12pt]

\begin{document}
\centerline{\bf A CONSTRUCTIVE DEMONSTRATION}

\centerline{\bf OF THE UNIQUENESS OF THE CHRISTOFFEL SYMBOL}

\bigskip\bigskip\bigskip

\centerline{{\bf Victor Tapia}\footnote{\tt
TAPIENS@CIENCIAS.CIENCIAS.UNAL.EDU.CO}}

\bigskip

\centerline{\it Departamento de F{\'\i}sica}

\centerline{\it Universidad Nacional de Colombia}

\centerline{\it Bogot\'a, Colombia}

\bigskip\bigskip\bigskip

The Christoffel symbol

\begin{equation}
\left\{{}^{\,\lambda}_{\mu\nu}\right\}({\bf
g})={1\over2}\,g^{\lambda\rho}\,(\partial_\mu g_{
\nu\rho}+\partial_\nu g_{\mu\rho}-\partial_\rho g_{\mu\nu})\,,
\label{1}
\end{equation}

\noindent is usually introduced (Spivak, 1975) as the connection which
solves the metricity condition

\begin{equation}
\nabla_\lambda g_{\mu\nu}=0\,.
\label{2}
\end{equation}

\noindent On the other hand, whether the Christoffel symbol is the only
connection which can be constructed from a symmetric second--rank tensor
$g_{\mu\nu}$ remains an open question. In this note we exhibit a
constructive demonstration of the uniqueness of the Christoffel symbol.

Let us start by reviewing some simple results of tensor calculus. Let
${\cal M}$ be an $n$--dimensional differentiable manifold. Several
geometric objects can be introduced over conveniently defined fibered
bundles based on ${\cal M}$. In order to classify them we adopt a
taxonomic approach, {\it cf.} (Visconti, 1992): a {\it tensor} is an
object which transforms like a {\it tensor}, etc. The previous definition
makes reference to the way in which an object transforms under changes of
local coordinates. Let $x^\mu$, $\mu=1,2,\cdots,n$, and $y^\alpha$,
$\alpha=1,2,\cdots,n$, be local coordinates on ${\cal M}$. Both sets are
related by $y^\alpha=y^\alpha(x^\mu)$, and the differential form of this
relation is

\begin{equation}
dy^\alpha=\left({{\partial y^\alpha}\over{\partial x^\mu}}\right)\,dx^\mu\,.
\label{3}
\end{equation}

\noindent Due to the intrinsic function theorem, this relation tells us
that $y^\alpha$ are functions of $x^\mu$, $y^\alpha=y^\alpha(x^\mu)$. In
order to express $x^\mu$ as functions of $y^\alpha$,
$x^\mu=x^\mu(y^\alpha)$, we need to invert eq. (\ref{3}). This can be
achieved if

\begin{equation}
\det\left({{\partial y^\alpha}\over{\partial x^\mu}}\right)\not=0\,.
\label{4}
\end{equation}

\noindent If this is the case we have an inverse matrix $(\partial
x/\partial y)$ and the relation (\ref{3}) can be inverted to

\begin{equation}
dx^\mu=\left({{\partial x^\mu}\over{\partial y^\alpha}}\right)\,dy^\alpha\,.
\label{5}
\end{equation}

\noindent Then $x^\mu=x^\mu(y^\alpha)$.

Now we are ready to proceed to our taxonomic classification of
geometrical objects. Relation (\ref{3}) provides the first example. Let
us consider a set of functions $v^\mu(x)$ and $v^\alpha(y)$ which are
related in the way given by (\ref{3})

\begin{equation}
v^\alpha(y)=\left({{\partial y^\alpha}\over{\partial
x^\mu}}\right)\,v^\mu(x)\,.
\label{6}
\end{equation}

\noindent This defines a contravariant vector. Accordingly, a set of
functions $v_\mu(x)$ and $v_\alpha(y)$ which are related by

\begin{equation}
v_\alpha(y)=\left({{\partial x^\mu}\over{\partial
y^\alpha}}\right)\,v_\mu(x)\,,
\label{7}
\end{equation}

\noindent defines a covariant vector. Let us now consider the functions
$\phi(x)=v_\mu(x)v^\mu(x)$ and $\phi(y)=v_\alpha(y)v^\alpha(y)$. From the
relations above we easily obtain

\begin{equation}
\phi(y)=\phi(x)\,.
\label{8}
\end{equation}

\noindent This relation defines a scalar. The relations above can be
extended to tensors of different rank and covariance. For example, a
covariant second--rank tensor $g_{\mu\nu}(x)$ is an object which
transforms as

\begin{equation}
g_{\alpha\beta}(y)=\left({{\partial x^\mu}\over{\partial
y^\alpha}}\right)\,\left({{\partial x ^\nu}\over{\partial
y^\beta}}\right)\,g_{\mu\nu}(x)\,.
\label{9}
\end{equation}

Let us now consider the derivative of the scalar function (\ref{8}). We
obtain

\begin{equation}
{{\partial\phi}\over{\partial y^\alpha}}(y)=\left({{\partial
x^\mu}\over{\partial y^\alpha}}
\right)\,{{\partial\phi}\over{\partial x^\mu}}(x)\,.
\label{10}
\end{equation}

\noindent Therefore, the ordinary derivative of a scalar function is a
vector. Next, let us consider the ordinary derivative of a covariant
vector. From (\ref{7}) we obtain

\begin{equation}
{{\partial v_\beta}\over{\partial y^\alpha}}(y)=\left({{\partial
x^\mu}\over{\partial y^\alpha}}
\right)\,\left({{\partial x^\nu}\over{\partial
y^\beta}}\right)\,{{\partial v_\nu}\over{\partial
x^\mu}}(x)+\left({{\partial^2x^\mu}\over{\partial y^\alpha\partial
y^\beta}}\right)\,v_\mu(x)\,.
\label{11}
\end{equation}

\noindent Therefore, this quantity is not a vector. Let us observe that
covariance is broken by a term linear in $v_\mu$. Let us therefore
introduce a new derivative

\begin{equation}
\nabla_\mu v_\nu=\partial_\mu v_\nu-{\Gamma^\lambda}_{\mu\nu}\,v_\lambda\,,
\label{12}
\end{equation}

\noindent where a term linear in $v_\mu$ is introduced to compensate the
wrong behaviour of the last term in (\ref{11}). Let us now impose that
this quantity be a tensor. The result is that ${\bf\Gamma}$ must
transform according to

\begin{equation}
{\Gamma^\gamma}_{\alpha\beta}(y)=\left({{\partial x^\mu}\over{\partial
y^\alpha}}\right)\, \left({{\partial x^\nu}\over{\partial
y^\beta}}\right)\,\left[\left({{\partial y^\gamma}\over{
\partial x^\lambda}}\right)\,{\Gamma^\lambda}_{\mu\nu}(x)-\left({{
\partial^2y^\gamma}\over{
\partial x^\mu\partial x^\nu}}\right)\right]\,.
\label{13}
\end{equation}

\noindent It is clear that ${\bf\Gamma}$ is not a tensor; this is not
unexpected since in order to compensate the non--tensor character of
$\partial v$ we need a non tensor object. $\nabla_\mu v_\nu$ is the
covariant derivative, and ${\bf\Gamma}$ is the connection. 

In Riemannian geometry the natural object is the metric tensor
$g_{\mu\nu}$. The metric is not related to the connection. However, we
may relate them through a metricity condition. The simplest metricity
condition is

\begin{equation}
\nabla^\Gamma_\lambda g_{\mu\nu}=0\,.
\label{14}
\end{equation}

\noindent The relation (\ref{14}) has the unique solution

\begin{equation}
{\Gamma^\lambda}_{\mu\nu}=\left\{{}^{\,\lambda}_{\mu\nu}\right\}({\bf
g})={1\over2}\,g^{
\lambda\rho}\,(\partial_\mu g_{\nu\rho}+\partial_\nu
g_{\mu\rho}-\partial_\rho g_{\mu\nu})\,,
\label{15}
\end{equation}

\noindent which is known as the Christoffel symbol. 

The Christoffel symbol is the connection which solves the metricity
condition (\ref{14}). It happens to be constructed from a second--rank
tensor. However, there is a further question we may ask: which is the
more general connection we can construct in terms of a second--rank
tensor? We must therefore look for a more general scheme providing a most
complete answer. 

Let us therefore consider the problem of how to construct a connection
${\Gamma^\lambda}_{\mu\nu}$ only from a second--rank tensor $g_{\mu\nu}$.
For this purpose let us remind that the connection is related to a
derivation and therefore has dimensions corresponding to an inverse
lenght, {\it i.e.}

\begin{equation}
{\rm dim}[{\bf\Gamma}]=L^{-1}\,.
\label{16}
\end{equation}

\noindent On the other hand, the derivation is a linear operation.
Therefore, the needed connection must be a linear combination of
derivatives of $g_{\mu\nu}$, namely,

\begin{equation}
{\Gamma^\lambda}_{\mu\nu}({\bf
g})=G^{\lambda\sigma\tau\rho}_{\mu\nu}({\bf g},\,{\bf g}^{-1})
\,\partial_\sigma g_{\tau\rho}\,,
\label{17}
\end{equation}

\noindent where ${\bf G}$ is a function depending only on $g_{\mu\nu}$
and its inverse $g^{\mu\nu}$, and therefore it is a tensor. The
transformation rule for ${\bf\Gamma}$ in (\ref{17}) is

\begin{eqnarray}
{\Gamma^\gamma}_{\alpha\beta}(y)&=&G^{\gamma\delta\epsilon\phi}_{\alpha\beta
}({\bf g},\,{\bf g}
^{-1})\,\partial_\delta g_{\epsilon\phi}\nonumber\\
&=&G^{\gamma\delta\epsilon\phi}_{\alpha\beta}({\bf g},\,{\bf
g}^{-1})\,\partial_\delta \left[
\left({{\partial x^\mu}\over{\partial
y^\epsilon}}\right)\,\left({{\partial x^\nu}\over{\partial
y^\phi}}\right)\,g_{\mu\nu}\right]\nonumber\\
&=&\left({{\partial y^\gamma}\over{\partial
x^\lambda}}\right)\,\left({{\partial y^\alpha}\over{
\partial x^\mu}}\right)\,\left({{\partial y^\beta}\over{\partial
x^\nu}}\right)\,{\Gamma^\lambda
}_{\mu\nu}(x)\nonumber\\
&&+2\,G^{\gamma\delta\epsilon\phi}_{\alpha\beta}({\bf g},\,{\bf
g}^{-1})\,\left({{\partial^2x^
\mu}\over{\partial y^\delta\partial y^\epsilon}}\right)\,\left({{\partial
x^\nu}\over{\partial y^\phi}}\right)\,g_{\mu\nu}\nonumber\\
&=&\left({{\partial y^\gamma}\over{\partial
x^\lambda}}\right)\,\left({{\partial y^\alpha}\over{
\partial x^\mu}}\right)\,\left({{\partial y^\beta}\over{\partial
x^\nu}}\right)\,{\Gamma^\lambda
}_{\mu\nu}(x)\nonumber\\
&&+2\,G^{\gamma\delta\epsilon\phi}_{\alpha\beta}({\bf g},\,{\bf
g}^{-1})\,g_{\phi\eta}\,\left(
{{\partial y^\eta}\over{\partial
x^\lambda}}\right)\,\left({{\partial^2x^\lambda}\over{\partial
y^\delta\partial y^\epsilon}}\right)\,.
\label{18}
\end{eqnarray}

\noindent Now we require that this object transform like a connection and
if we compare with (\ref{13}) we obtain

\begin{equation}
2\,G^{\gamma\delta\epsilon\phi}_{\alpha\beta}({\bf g},\,{\bf
g}^{-1})\,g_{\phi\eta}\,\left({{
\partial y^\eta}\over{ \partial
x^\lambda}}\right)\,\left({{\partial^2x^\lambda}\over{\partial y
^\delta\partial y^ \epsilon}}\right)=\left({{\partial
y^\gamma}\over{\partial x^\lambda}}\right)
\,\left({{\partial^2x^\lambda}\over{\partial y^\alpha\partial
y^\beta}}\right)\,.
\label{19}
\end{equation}

\noindent Therefore

\begin{equation}
\left[4\,G^{\gamma\delta\epsilon\phi}_{\alpha\beta}({\bf g},\,{\bf
g}^{-1})\,g_{\phi\eta}-\delta
^\gamma_\eta\,\left(\delta^\delta_\alpha\,\delta^\epsilon_\beta+\delta^
\delta_\beta\,\delta^
\epsilon_\alpha\right)\right]\,\left({{\partial y^\eta}\over{\partial
x^\lambda}}\right)\,\left(
{{\partial^2x^\lambda}\over{\partial y^\delta\partial
y^\epsilon}}\right)=0\,.
\label{20}
\end{equation}

\noindent What must be zero is the symmetric part of the square bracket
and then we obtain

\begin{equation}
2\,G^{\gamma\delta\epsilon\phi}_{\alpha\beta}({\bf g},\,{\bf
g}^{-1})\,g_{\phi\eta}+2\,G^{\gamma
\epsilon\delta\phi}_{\alpha\beta}({\bf g},\,{\bf
g}^{-1})\,g_{\phi\eta}-\delta^\gamma_\eta\,
\left(\delta^\delta_\alpha\,\delta^\epsilon_\beta+\delta^\delta_\beta\,\delta
^\epsilon_\alpha
\right)=0\,.
\label{21}
\end{equation}

\noindent Considering cyclic permutations of the indices
$\delta\epsilon\phi$ we arrive to a set of equations which allow to
determine ${\bf G}$. The solution is

\begin{eqnarray}
&&G^{\gamma\delta\epsilon\phi}_{\alpha\beta}({\bf g},\,{\bf
g}^{-1})\nonumber\\
&=&{1\over4}\,\left[g^{\gamma\phi}\,\left(\delta^\delta_\alpha\,\delta^
\epsilon_\beta+\delta^
\delta_\beta\,\delta^\epsilon_\alpha\right)+g^{\gamma\delta}\,\left(\delta^
\epsilon_\alpha\,
\delta^\phi_\beta+\delta^\phi_\beta\,\delta^\epsilon_\alpha\right)-g^{\gamma
\epsilon}\,\left(
\delta^\delta_\alpha\,\delta^\phi_\beta+\delta^\phi_\beta\,
\delta^\delta_\alpha\right)\right]
\,.\nonumber\\
&&\label{22}
\end{eqnarray}

\noindent If we now replace this result in (\ref{16}) we obtain that the
connection is precisely the Christoffel symbol, namely eq. (\ref{15}). 

One can inmediately check that

\begin{equation}
\nabla^g_\lambda g_{\mu\nu}\equiv0\,,
\label{23}
\end{equation}

\noindent where $\nabla^g$ is the covariant derivative constructed with
the Christoffel symbol. 

We have shown that the Christoffel symbol is the only connection which
can be constructed from a symmetric second--rank tensor. The metricity
condition (\ref{23}) appears just as a side result and not as the
starting point for the construction of the Christoffel symbol.

\bigskip\bigskip

\centerline{\bf References}

\begin{enumerate}

\item M. Spivak, {\it A comprehensive introduction to differential
geometry} (Publish or Perish, Boston, 1975).

\item A. Visconti, {\it Introductory differential geometry for
physicists} (World Scientific, River Edge, 1992).

\end{enumerate}

\end{document}